\documentclass[twocolumn]{autart}
\usepackage{amsfonts,amssymb}
\usepackage{amsmath, natbib}
\begin{document}

\begin{frontmatter}
\title{Further Results on Strict Lyapunov Functions for\\ Rapidly
Time-Varying Nonlinear Systems\thanksref{footnoteinfo}}

\thanks[footnoteinfo]{Corresponding Author: Fr\'ed\'eric Mazenc. The second author was  supported by NSF/DMS
Grant 0424011.  }

\author[1]{Fr\'ed\'eric Mazenc}\ead{Frederic.Mazenc@ensam.inra.fr},
\author[2]{Michael Malisoff}\ead{malisoff@lsu.edu},
\author[3]{Marcio S. de Queiroz}\ead{dequeiroz@me.lsu.edu}

\address[1]{Projet MERE INRIA-INRA, UMR Analyse des Syst\`emes et
Biom\'etrie INRA, 2 pl. Viala, 34060 Montpellier, France}
\address[2]{Department of
Mathematics, Louisiana State University, Baton Rouge, LA 70803-4918}
\address[3]{Department of Mechanical Engineering, Louisiana
State University, Baton Rouge, LA 70803-6413}

\begin{keyword}time-varying systems, input-to-state stability,
Lyapunov function constructions\end{keyword}

\begin{abstract}
We explicitly construct  global strict Lyapunov functions for
rapidly time-varying nonlinear control  systems.  The Lyapunov
functions we construct are expressed  in terms of oftentimes more
readily available Lyapunov functions for the limiting dynamics which
we assume are uniformly globally asymptotically stable. This leads
to new sufficient conditions for uniform global exponential, uniform
global asymptotic, and input-to-state stability of fast time-varying
dynamics.  We also construct strict Lyapunov functions for our
systems using a strictification approach. We illustrate our results
using several examples.
\end{abstract}
\end{frontmatter}

\section{Introduction}
The stabilization of nonlinear and nonautonomous control systems,
and the construction of their Lyapunov functions,  are challenging
problems that are of significant ongoing  interest in the context of
robustness analysis and controller design; see
\cite{MM05} and \cite{MB04}. One popular approach to guaranteeing
stability of nonautonomous systems is the {\em averaging method} in
which exponential stability of an appropriate {\em autonomous}
system implies exponential stability of the original dynamics when
its time variation  is sufficiently fast. See \cite{K02} for related
results.

The preceding results were extended to more general rapidly
time-varying systems of the form
\begin{equation}
\label{m1} \dot{x} \; = \; f(x,t,\alpha t),\; \; x\in {\mathbb
R}^n,\; t\in {\mathbb R},\; \alpha>0
\end{equation}
in \cite{PA02}, where uniform (local) exponential stability of
(\ref{m1}) was proven for large values of the constant $\alpha>0$,
assuming a suitable limiting dynamics
\begin{equation}
\label{m2} \dot{x} \; = \; \bar{f}(x,t)
\end{equation}
for (\ref{m1}) is uniformly exponentially stable.  (We specify the
choice of $\bar f$ in  our main theorem and  examples below.) This
generalized a result from \cite[pp. 190-5]{H80} on  a class of
systems (\ref{m1}) satisfying certain periodicity or almost
periodicity conditions. The main arguments of \cite{PA02} use
(partial) averaging  but do not lead to explicit Lyapunov functions
for (\ref{m1}).

In this work, we pursue a very different approach. Instead of
averaging, we
 explicitly construct a
family of Lyapunov functions for (\ref{m1}) {}in terms of more
readily available Lyapunov functions for (\ref{m2}), which we again
assume is asymptotically stable. {}While \cite{PA02} assumes
(\ref{m2}) is uniformly  {}locally exponentially stable, we consider
the case where (\ref{m2}) is assumed uniformly globally
asymptotically stable (UGAS), in which case our conclusion is that
(\ref{m1}) is UGAS (but not necessarily exponentially stable) when
$\alpha>0$ is sufficiently large. While global exponential and
global asymptotic stabilities are equivalent for {\em autonomous}
systems under a coordinate change in certain dimensions, the
coordinate changes  are not explicit and so do not lend themselves
to explicit Lyapunov function constructions;  see \cite{GSW99}.  The
importance of the problem we consider lies in the ubiquity of
rapidly time-varying systems in practical applications (e.g.,
suspended pendulums subject to vertical vibrations of small
amplitude and high {}frequency, Raleigh's and Duffing's equations
from \cite[Chatper 10]{K02}, {}and systems arising in identification
{}in \cite{PA02} and below) and the essentialness of Lyapunov
functions in robustness analysis and controller design for these
applications.

In particular, we show that assumptions similar to those of
\cite[Theorem 3]{PA02}  imply that (\ref{m1}) is uniformly globally
 {}exponentially stable;  our Lyapunov
function constructions are new even in this particular case and our
results are complementary to those of \cite{PA02}. The Lyapunov
functions  we construct are also input-to-state stable (ISS) or
integral ISS Lyapunov functions for
\begin{equation}
\label{fastcontrol} \dot{x}\; =\; f(x,t,\alpha t)+g(x,t,\alpha
t)u\end{equation} under  appropriate  conditions on $f$ and $g$; see
Remark \ref{gissue}.

  In Section \ref{dl}, we provide the relevant
definitions and lemmas. In Section \ref{thm}, we present our main
sufficient conditions for uniform  global asymptotic and exponential
stability of (\ref{m1}), and for the stability of
(\ref{fastcontrol}), in terms of limiting dynamics (\ref{m2}). This
theorem leads to explicit constructions of Lyapunov functions for
(\ref{m1}) and (\ref{fastcontrol})  in terms of Lyapunov functions
for (\ref{m2}), in Theorem \ref{cor1}.  In Section \ref{altthmsec},
we provide an alternative Lyapunov function construction theorem for
(\ref{m1}) not involving any limiting dynamics.  We prove our main
results in Sections \ref{proof} and \ref{altthm}. We illustrate our
theorems in Sections \ref{illu} and \ref{altillus}
  using a
friction control model and other examples.
 We close in Section \ref{concl} by
summarizing our findings.

\section{Assumptions, definitions, and lemmas}
\label{dl} We study  (\ref{m1}) (which includes dynamics (\ref{m2})
with no $\alpha$ dependence, as special cases) in which  we always
assume $f$ is  continuous in time $t\in {\mathbb R}:=(-\infty,
+\infty)$, continuously differentiable ($C^1$) in  $x\in {\mathbb
R}^n$,
  null at $x=0$ meaning
\begin{equation}\label{nullness}f(0,t,\alpha t)= \bar f(0,t)=
0\; \; \; \forall t\in {\mathbb R}, \; \alpha>0\end{equation} and
{\em forward complete}, i.e., for each $\alpha>0$, $x_o\in {\mathbb
R}^n$, and $t_o\in {\mathbb R}_{\ge 0}:=[0,\infty)$   there {}is a
unique trajectory $[t_0,\infty)\ni t\mapsto \phi(t; t_o, x_o)$ for
(\ref{m1}) (depending  in general on the constant $\alpha>0$) that
satisfies $x(t_o)=x_o$. We set ${\mathbb N}=\{1,2,3,\ldots\}$ and
let ${\mathbb Z}$ denote the set of all integers.  We say
$N:{\mathbb R}_{\scriptstyle\ge 0}\to {\mathbb R}_{\scriptstyle\ge
0}$ is {\em of class ${\mathcal  M}$} and write $N\in {\mathcal  M}$
provided
\begin{equation}
\label{m3} \lim_{\eta\to +\infty}\eta N(\eta)=0. \end{equation} A
continuous function $\delta:{\mathbb R}_{\scriptstyle\ge 0}\to
{\mathbb R}_{\scriptstyle\ge 0}$ is {\em positive definite} provided
it is zero only at zero.  A positive definite  function $\delta$ is
{\em of class ${\mathcal  K}$} (written $\delta\in {\mathcal  K}$)
provided it is strictly increasing;  if in addition $\delta$ is
unbounded, then we say that $\delta$ is {\em of class ${\mathcal
K}_\infty$} and write $\delta\in {\mathcal K}_\infty$. A continuous
function  $\beta:{\mathbb R}_{\scriptstyle\ge 0}\times {\mathbb
R}_{\scriptstyle\ge 0}\to {\mathbb R}_{\scriptstyle\ge 0}$ is {\em
of class $\mathcal{KL}$} (written $\beta\in \mathcal{KL}$) provided
(a) $\beta(\cdot, t)\in {\mathcal K}_\infty$ for all $t\ge 0$, (b)
$\beta(s,\cdot)$ is nonincreasing for all $s\ge 0$, and (c) for each
$s\ge 0$, $\beta(s,t)\to 0$ as $t\to +\infty$.  A positive definite
function $\delta$  is called $o(s)$ provided $\delta(s)/s\to 0$ as
$s\to +\infty$.  {}We always assume there exists $\rho\in
\mathcal{K}_\infty$ such that $|f(x,t,\alpha t)|\le \rho(|x|)$
everywhere.

We next define our stability properties for (\ref{m2}).  The same
definitions apply for  (\ref{m1}) for any choice of the constant
$\alpha>0$. We call  (\ref{m2})  {\em uniformly globally
asymptotically stable (UGAS)} provided there exists $\beta\in
\mathcal{KL}$  such that
\begin{equation}
\label{gasdef} \begin{array}{l}\! \! \! \! \! |\phi(t; t_o, x_o)|\,
\le\, \beta(|x_o|, t-t_o)\;   \forall t\ge t_o\ge 0, \, x_o\in
{\mathbb R}^n\end{array}\end{equation} where $\vert\cdot\vert$ is
the usual Euclidean norm and $\phi$ is the flow map for (\ref{m2}).
We call (\ref{m2}) {\em uniformly globally exponentially stable
(UGES)} provided there exist constants $D>1$ and $\lambda>0$ such
that (\ref{gasdef}) is satisfied with the choice
\begin{equation}
\label{betachoice} \beta(s,t)\, =\, D s e^{-\lambda t}.
\end{equation}
The converse Lyapunov function theorem  implies  (\ref{m2}) is UGAS
if and only if  it has a {\em(strict) Lyapunov function}, i.e., a
$C^1$  $V:{\mathbb R}^n\times {\mathbb R}_{\scriptstyle\ge 0}\to
{\mathbb R}_{\scriptstyle\ge 0}$ that admits $\delta_1,\delta_2\in
{\mathcal  K}_\infty$ and $\delta_3\in {\mathcal K}$ such that for
all $t\in {\mathbb R}_{\scriptstyle\ge 0}$ and $\xi\in {\mathbb
R}^n$, we have both (L1) $\delta_1(|\xi|)\le V(\xi,t) \le
\delta_2(|\xi|)$ and (L2) $V_t(\xi,t) + V_\xi(\xi,t) \,
\bar{f}(\xi,t)  \leq  -\delta_3(|\xi|)$, where the subscripts on $V$
{}indicate partial gradients;  see \cite{BR05}. When (\ref{m2}) is
UGES,  the proof of \cite[Theorem 4.14]{K02} shows:
\begin{lem}
\label{k02} Assume (\ref{m2}) satisfies the UGES condition
(\ref{gasdef})-(\ref{betachoice}) for some constants $D>1$ and
$\lambda>0$ and that there exists $K>\lambda$ such that $|(\partial
\bar f/\partial \xi)(\xi,t)|\le K$ for all $\xi\in {\mathbb R}^n$
and $t\in {\mathbb R}_{\scriptstyle\ge 0}$.  Then (\ref{m2}) admits
a Lyapunov function $V$ and constants $c_1,c_2,c_3>0$ such that
\begin{equation}
\label{m5}\renewcommand{\arraystretch}{1.25}
\begin{array}{l}
 c_1|\xi|^2 \; \leq \; V(\xi,t) \; \leq \; c_2|\xi|^2\; , \; \;
\left|V_\xi(\xi,t)\right| \; \leq \; c_3 |\xi|\; ,\\ {\rm and}\; \;
\; V_t(\xi,t) + V_\xi(\xi,t) \bar{f}(\xi,t) \; \leq \; -
|\xi|^2\end{array}
\end{equation}
hold for all $t\in {\mathbb R}_{\scriptstyle\ge 0}$ and {}$\xi\in
{\mathbb R}^n$.
\end{lem}

\begin{defn}
\label{compatt} Given $\delta\in {\mathcal  K}$, the dynamics
(\ref{m2}) is said to be {\em $\delta$-compatible} provided it
admits a Lyapunov function $V\in C^1$ and two constants $\bar c\in
(0,1)$, $\bar{\bar c}
>0$ such that:
\begin{itemize}\item\begin{itemize}\addtolength{\itemsep}{.25\baselineskip}
\item[$P_1$]\ \  $V_t(\xi,t) + V_\xi(\xi,t) \, \bar{f}(\xi,t) \;
\leq \; -\bar c\, \delta^2(|\xi|)\; \; \; \; \forall \xi,t$.
\item[$P_2$]\ \  $|V_\xi(\xi,t)| \; \leq \; \delta(|\xi|)\;
 {\rm and}\;  |\bar{f}(\xi,t)|\; \le\;  \delta(|\xi|/2)\; \;
\; \;  \forall \xi,t$. \item[$P_3$]\ \  $\delta(s)\; \le\; \bar{\bar
c}\, s\; \; \; \;  \forall s\ge 0$.\end{itemize}
\end{itemize}\smallskip\end{defn}
\begin{rem}\rm\label{compat}Note the {\em asymmetry} in the bounds on $|V_{\scriptscriptstyle \xi}|$ and
$|\bar f|$ in $P_2$. If (\ref{m2}) satisfies the assumptions of
Lemma \ref{k02}, then it is $\delta$-compatible with
$\delta(s)=(c_3+2K)s$. However, by varying $\delta$ (including cases
where $\delta$ is bounded), {}one also finds a rich class of
non-UGES $\delta$-compatible dynamics; see e.g. Section \ref{exa1}
below.\end{rem}

We also consider the nonautonomous {\em control system}
\begin{equation}
\label{nonauto} \dot x=F(x,t,u)
\end{equation}
which we always assume is continuous in all variables and $C^1$ in
$x$ with $F(0,t,0)\equiv 0$, and whose solution for  a given control
function ${\mathbf u}\in{\mathcal  U}(:=$all measurable locally
essentially bounded functions
 $[0,\infty)\to {\mathbb R}^m$) and  given initial condition
 $x(t_o)=x_o$ we denote by $t\mapsto \phi(t; t_o,x_o,{\mathbf u})$.
 We
 always assume (\ref{nonauto}) is {\em forward complete},
 i.e.,
 all trajectories $\phi(\cdot; t_o, x_o, {\mathbf u})$  so defined have domain $[t_o,+\infty)$.
 We
 next recall the input-to-state stable (ISS) and integral
 input-to-state stable
 (iISS) properties from
 \cite{S89} and \cite{S98}.  Let $|{\mathbf
 u}|_I$ denote the essential supremum of  ${\mathbf u}\in {\mathcal  U}$
restricted to any interval $I\subseteq {\mathbb R}_{\scriptstyle\ge
0}$.
\begin{defn}
\label{iss} (a) We say that (\ref{nonauto}) is {\em ISS} provided
there exist $\gamma\in \mathcal{K}_\infty$ and $\beta\in
\mathcal{KL}$ for which
\begin{equation}
\label{issdef} |\phi(t; t_o, x_o, {\mathbf u})|\; \le\; \beta(|x_o|,
t-t_o)+\gamma\left(|{\mathbf u}|_{[t_o, t+t_o]}\right)
\end{equation} holds when $t\ge t_o\ge 0$, $x_o\in {\mathbb
R}^n$, and ${\mathbf u}\in {\mathcal  U}$. If in addition $\beta$
has the form (\ref{betachoice}), then we say that (\ref{nonauto}) is
{\em input-to-state exponentially stable (ISES)}.  (b) We say that
(\ref{nonauto}) is {\em iISS} provided there exist $\mu, \gamma\in
{\mathcal  K}_\infty$ and $\beta\in \mathcal{KL}$ such that
\[\mu(|\phi(t; t_o, x_o, {\mathbf u})|)\; \le\;  \beta(|x_o|,
t-t_o)+\int_{t_o}^{t_o+t}\gamma(|{\mathbf u}(s)|)\, {\rm d}s
\]
 holds when $t\ge t_o\ge 0$, $x_o\in {\mathbb R}^n$, and ${\mathbf u}\in {\mathcal
U}$.
\end{defn}

A function $V: {\mathbb R}^n\times {\mathbb R}_{\ge 0}\to {\mathbb
R}_{\ge 0}$ is called {\em uniformly positive definite} provided
$s\mapsto \inf\{V(x,t):t\ge 0, |x|=s\}$ is positive definite in
which case we write $V\in \text{UPD}$.  We call $V\in \text{UPD}$
{\em uniformly proper and positive definite} provided there exist
$\delta_1,\delta_2\in {\mathcal K}_\infty$ such that condition (L1)
above is satisfied in which case we write $V\in\text{UPPD}$. The
following Lyapunov function notions  agree with the  usual ISS and
iISS Lyapunov function definitions when (\ref{nonauto}) is
autonomous, because functions $\chi\in {\mathcal K}_\infty$ are
invertible. Note that $\nu$ in (\ref{iissd})    {\em need not} be of
class ${\mathcal  K}$.

\begin{defn}
\label{corr} Let $V\in C^1\cap {\rm UPPD}$. (a) We call $V$ an {\em
ISS Lyapunov function} for (\ref{nonauto}) provided there exist
$\chi,\delta_3\in {\mathcal K}_\infty$ such that for all $t\in
{\mathbb R}_{\scriptstyle\ge 0}$, $\xi\in {\mathbb R}^n$, and $u\in
{\mathbb R}^m$:
$[|u|\le \chi(|\xi|)  \Rightarrow   V_t(\xi,t) + V_\xi(\xi,t) \,
F(\xi,t,u)  \leq  -\delta_3(|\xi|)]$.
%
(b) We call $V$ an {\em iISS Lyapunov function} for (\ref{nonauto})
provided there exist $\Delta\in {\mathcal K}_\infty$ and a positive
definite function $\nu:{\mathbb R}_{\scriptstyle\ge 0}\to{\mathbb
R}_{\scriptstyle\ge 0}$   such that
\begin{equation}\label{iissd}V_t(\xi,t) + V_\xi(\xi,t) \, F(\xi,t,u) \; \leq \;
-\nu(|\xi|)+\Delta(|u|)
\end{equation}
holds for all $t\in {\mathbb R}_{\scriptstyle\ge 0}$, $\xi\in
{\mathbb R}^n$, and $u\in {\mathbb R}^m$.
\end{defn}
Since (\ref{nonauto}) has an ISS Lyapunov function when it is ISS
(by the arguments of \cite{SW95}), the proof of
 \cite[Theorem 1]{ASW00} shows that if (\ref{nonauto}) is ISS {}and autonomous, then it is also iISS, but not
conversely, since e.g. $\dot x=-\arctan(x)+u$ is iISS but not ISS.
The next lemma follows from the arguments used in \cite{ASW00},
\cite{ELW00}, and \cite{S89}.
\begin{lem}
\label{suff} If (\ref{nonauto}) admits an ISS (resp., iISS) Lyapunov
function, then it is ISS (resp., iISS).\end{lem}

\section{Statements and discussions of main results}

\subsection{Main theorem and Lyapunov function construction} \label{thm} We show that {}hypotheses
similar to those of \cite[Theorem 3]{PA02}  ensure  that (\ref{m1})
is UGES. In fact, we show {}our conditions imply
\begin{equation}
\label{dis} \dot{x} \; = \; f(x,t,\alpha t)+u,\; \; \; x\in {\mathbb
R}^n, \; u\in {\mathbb R}^n
\end{equation}
is ISS when $\alpha>0$ is sufficiently large; see Remark
\ref{gissue} below for results on the more general systems
(\ref{fastcontrol}). Our main assumption  will be: There exist
$\delta\in {\mathcal  K}$, a $\delta$-compatible dynamics
(\ref{m2}), and $N\in {\mathcal  M}$ (cf.  (\ref{m3}) above) such
that for all $x\in {\mathbb R}^n$, all $r\in {\mathbb R}$ and
sufficiently large $\eta>0$,
\begin{equation}\label{relate}   \left|\int_{r - \frac{1}{\eta}}^{r +
\frac{1}{\eta}} \left\{f(x,l,\eta^2 l) - \bar{f}(x,l)\right\} {\rm
d}l \right|   \leq   \delta(|x|/2)  N(\eta)
\end{equation}
{}which is similar to  \cite[Property 2]{PA02}  in the special case
where $\delta(s)=2s$. Two more advantages of our result are  (a) it
applies to cases where  (\ref{m2}) is UGAS but not necessarily UGES
(cf. Section \ref{exa1} below) and (b) its proof leads to explicit
 Lyapunov functions for {}(\ref{dis}) (cf. Theorem \ref{cor1} below).
See also  Theorem \ref{altthm} below for cases where $\partial
f/\partial x$ is not necessarily globally bounded.  Recall the
definition of compatibility (in Definition \ref{compatt}) and the
requirement $N\in \mathcal{M}$  in (\ref{m3}).
\begin{thm}
\label{lef2} Consider a system (\ref{m1}). Assume there exist
$\delta\in {\mathcal  K}$, a $\delta$-compatible UGAS system
(\ref{m2}), two constants $\eta_o>0$ and $K>1$,  and $N\in {\mathcal
M}$ such that (\ref{relate}) holds whenever $\eta \geq \eta_o$,
$x\in {\mathbb R}^n$ and $r\in {\mathbb R}$ and such that:
\begin{equation} \label{m16}\renewcommand{\arraystretch}{1.5}\begin{array}{l} \left|\frac{\partial
\bar{f}}{\partial x}(x,t)\right| \, \leq \, K \; \; , \; \;
\left|\frac{\partial f}{\partial x}(x,t,\alpha t)\right| \, \leq \,
K\; \; ,\; \; {\rm and}\\ |f(x,t,\alpha t)|\, \le\,  \delta(|x|/2)\;
\; \; \; \forall t\in{\mathbb R}, \, x\in {\mathbb R}^n,\,
\alpha>0.\end{array}
\end{equation}
  Then there is  a
constant $\underline{\alpha}>0$ such that for all constants $\alpha
\geq \underline{\alpha}$, (\ref{m1}) is UGAS  and (\ref{dis}) is
iISS. If in addition $\delta\in {\mathcal K}_\infty$, then
(\ref{dis}) is ISS for all constants $\alpha \geq
\underline{\alpha}$.
 In the special case where (\ref{m2}) is UGES,
 (\ref{m1})
is  UGES for all constants $\alpha \geq \underline{\alpha}$ and
(\ref{dis}) is  ISES for all constants $\alpha \geq
\underline{\alpha}$.
\end{thm}

By (\ref{nullness}),  the condition $|f(x,t,\alpha t)| \le
\delta(|x|/2)$ in (\ref{m16}) is redundant when $\delta$ has the
form  $\delta(s)=\bar rs$ for a constant $\bar r>0$, since $\bar r$
can always be enlarged. Our proof of  Theorem \ref{lef2} in Section
\ref{proof} below  will also show:

\begin{thm}\label{cor1}

Let the hypotheses of Theorem  \ref{lef2} hold for some $\delta\in
{\mathcal  K}$, and $V\in C^1$ be {}a Lyapunov function for
(\ref{m2}) satisfying the requirements of Definition \ref{compatt}.
Then there exists a constant $\underline{\alpha}>0$ such that for
all constants $\alpha>\underline{\alpha}$,
\[\begin{array}{l}
V^{\scriptscriptstyle [\alpha]}(\xi,t):=\\
V\left(\xi-\frac{\sqrt{\alpha}}{2}\int_{
t-\frac{2}{\sqrt{\alpha}}}^t\int_s^t\{f(\xi,l,\alpha
l)-\bar{f}(\xi,l)\}\, {\rm d} l\, {\rm d}s,t\right)\end{array}\] is
a Lyapunov function for (\ref{m1}) and an iISS Lyapunov function for
(\ref{dis}).
 If also $\delta\in {\mathcal  K}_\infty$, then
$V^{\scriptscriptstyle [\alpha]}$ is also an ISS Lyapunov function
for
 (\ref{dis}) for all constants
$\alpha>\underline{\alpha}$.
\end{thm}

\subsection{Alternative result}
\label{altthmsec} \label{alt} The proof of Theorem \ref{lef2}
constructs strict Lyapunov functions for (\ref{m1}) in terms  of
strict Lyapunov functions for the limiting dynamics (\ref{m2}). It
is natural to inquire whether one can instead construct strict
Lyapunov functions for (\ref{m1}) by strictifying nonstrict Lyapunov
functions for (\ref{m1});  see \cite{MM05} where the strictification
approach was applied to nonautonomous systems that are not rapidly
time-varying. In this section, we extend this approach to cover
(\ref{m1}). Our strictification result has the advantages in certain
situations that (a) it does not require any knowledge of  limiting
dynamics, (b) it allows the derivative of
 the nonstrict Lyapunov function to be zero or even
positive at some points, and (c) it does not require (\ref{m1}) to
be globally Lipschitz in the state. In the rest of the section, we
focus on systems with no controls, but the extension to the control
system (\ref{fastcontrol}) for appropriate $g$ can be done by
similar arguments.
 We  assume:
\begin{itemize}\item[]\begin{itemize}
\item[H.]
There exist $V\in C^1\cap {\rm UPPD}$, $W\in {\rm UPD}$, a $C^1$
function $\Theta: {\mathbb R}^n\times {\mathbb R}_{\ge 0}\to
{\mathbb R}$, a bounded continuous function $p:{\mathbb
R}\to{\mathbb R}$, and constants $c,T>0$ such that  for all $x\in
{\mathbb R}^n$, $t\ge 0$, $\alpha>0$, and $k\in  {\mathbb Z}$, we
have:\smallskip\smallskip
\begin{itemize}
\item[H1.]
$V_t+V_xf(x,t,\alpha t) \leq - W(x,t) + p(\alpha t) \Theta(x,t)$
\item[H2.]  $\int_{kT}^{(k+1)T}p(r){\rm d}r=0$\smallskip
\item[H3.] $V\ge c|\Theta|$, $W\ge c\max\{|\Theta|,
|\Theta_t+\Theta_xf(x,t,\alpha t)|\}$
\end{itemize}
\end{itemize}\end{itemize}
Here and in the sequel, we  omit the argument  $(x,t)$ of $V$, $W$,
$\Theta$, and the partial gradients of $V$ and $\Theta$ whenever
this would not lead to  confusion.
In Section \ref{altproof}, we prove:
\begin{thm}
\label{altthm} If Assumption H. holds, then there exists a constant
$\underline \alpha>0$ such that for all constants $\alpha\ge
\underline \alpha$,
\begin{equation}
\label{jb5}
\begin{array}{l}
U^{\scriptscriptstyle [\alpha]}(x,t)\\ = V(x,t) - \left(\int_{t -
1}^{t} \left(\int_{s}^{t} p(\alpha l) dl\right) ds\right)
\Theta(x,t)\end{array}
\end{equation}
is a Lyapunov function for (\ref{m1}).  In particular, (\ref{m1}) is
UGAS for all constants $\alpha\ge \underline\alpha$.
\end{thm}

\section{Proof of  Theorems \ref{lef2} and \ref{cor1} and
remarks}\label{proof} To make our arguments easy to follow, we first
outline our method for proving these theorems.  First, we give the
proof of Theorem \ref{lef2} for the special case where (\ref{m2}) is
UGAS and $\delta\in \mathcal{K}_\infty$ which includes the proof of
Theorem \ref{cor1} for the  $\delta\in \mathcal{K}_\infty$ case.
Then we indicate the changes required if $\delta\in \mathcal{K}$ is
bounded.  Finally, we specialize to the special case where
$\delta(s)=\bar r s$ for some constant $\bar r>0$ which will prove
the ISES assertion of Theorem \ref{lef2}.

We first assume that (\ref{m2}) is UGAS and $\delta\in {\mathcal
K}_\infty$,  and we prove the ISS property for (\ref{dis}) for large
constants $\alpha>0$. In what follows, we assume all inequalities
and equalities hold {}globally, unless otherwise indicated. Let
$\eta_o$ be as in the statement of the theorem, and fix $\alpha =
\eta^2$ with $\eta \geq \eta_o$, ${\mathbf u}\in {\mathcal U}$, and
a trajectory $x(t)$ for (\ref{dis}) and ${\mathbf u}$, with
arbitrary initial condition.
   Set
\begin{equation}
\label{ze1} z(t) = x(t) + R_\alpha(x(t),t),
\end{equation}
where
\[
 R_\alpha(x,t)= - \frac{\eta}{2}\int_{t -
2/\eta}^{t}\int_{s}^{t} \left\{f(x,l,\eta^2 l) -
\bar{f}(x,l)\right\}  {\rm d} l\,  {\rm d}s.\] This is well-defined
since we are assuming our dynamics are forward complete. Set
{}\[p(t,l)=f(x(t),l,\eta^2 l) - \bar{f}(x(t),l).\] With $p$ so
defined, one easily checks (as was done e.g. in \cite{MM05}) that
for any $\tau>0$ and $t\in {\mathbb R}$,
\begin{equation}
\begin{array}{l}\label{fua}
 \displaystyle\frac{d}{dt}\int_{t-\tau}^t\int_s^tp(t,l)\, \, {\rm d} l\, {\rm
d}s=\\ \tau p(t,t)-\int_{t-\tau}^tp(t,l)\, \, {\rm d} l +\;
\displaystyle\int_{t-\tau}^t\int_s^t\frac{\partial p}{\partial
t}(t,l)\, \, {\rm d} l\, {\rm d}s\end{array}\end{equation}
and\begin{equation}
\begin{array}{l}\label{fub}
   \left\vert \int_{t-\tau}^t\int_s^tp(t,l)\, \, {\rm d}
l\, {\rm d}s\right\vert \; \le\;
\frac{\tau^2}{2}\displaystyle\max_{t-\tau\le l\le t} |p(t,l)|.
\end{array}\end{equation}
Taking $\tau=2/\eta$, (\ref{fua})  multiplied through by $-\eta/2$,
$f(x(t),t,\eta^2 t)-p(t,t)\equiv \bar f(x(t),t)$,  and (\ref{ze1})
give \begin{equation}\!  \begin{array}{l} \label{along} \dot{z}(t) =
\bar{f}(z(t),t) + \left(\bar{f}(x(t),t) - \bar{f}(z(t),t)\right) +
{\mathbf u}(t)\\  +\frac{\eta}{2}\int_{t - \frac{2}{\eta}}^{t}
p(t,l) \, {\rm d} l- \frac{\eta}{2}\left\{\int_{t -
\frac{2}{\eta}}^{t}\int_{s}^{t} \left(\frac{\partial f}{\partial
x}(x(t),l,\eta^2l)\right.\right.\\\left.\left. - \frac{\partial
\bar{f}}{\partial x}(x(t),l)\right) \, {\rm d} l \, {\rm
d}s\right\}\left(f(x(t), t,\eta^2t)+ {\mathbf
u}(t)\right).\end{array}\end{equation}

Let $V$, $\delta_1$, and $\delta_2$ satisfy the requirements
$P_1$-$P_3$ from our compatibility condition (in Definition
\ref{compatt}) and (L1).
 By $P_1$ with $\xi=z(t)$ {}and (\ref{along}),
 the derivative of $V(z,t)$ along the time-varying map $z(t)$ defined in (\ref{ze1}) (which we denote simply by $\dot V$ in the sequel) satisfies
\[\begin{array}{l}\renewcommand{\arraystretch}{1}
 \dot V  \; \le\;     -\bar c\, \delta^2(|z(t)|) + V_\xi(z(t),t)
\left(\bar{f}(x(t),t)-\bar{f}(z(t),t)\right)\nonumber\\
+ \frac{\eta}{2}V_\xi(z(t),t)\int_{t - 2/\eta}^{t} p(t, l) \, {\rm
d} l \nonumber
\\
 - \frac{\eta}{2}V_\xi(z(t),t)\\\times
\left[\int_{t - 2/\eta}^{t}\int_{s}^{t} \left(\frac{\partial
f}{\partial x}(x(t),l,\eta^2l) - \frac{\partial \bar{f}}{\partial
x}(x(t),l)\right) \, {\rm d} l\, {\rm d}s\right]\nonumber\\ \times
\left(f(x(t), t,\eta^2t)+ {\mathbf u}(t)\right)+V_\xi(z(t),t)\,
\mathbf{u}(t).\nonumber
\end{array}\]
 We deduce from (\ref{relate}), (\ref{m16}), (\ref{ze1}), and $P_2$  that
  \[\renewcommand{\arraystretch}{1.75}
\begin{array}{rcl}
\dot{V} & \leq & - \bar c\, \delta^2(|z(t)|) + K\delta(|z(t)|)|x(t)
- z(t)|
 \\&&+ \frac{\eta}{2}\, \delta(|z(t)|)\left|\int_{t - 2/\eta}^{t}
p(t,l)
 \, {\rm d} l
\right|+\delta(|z(t)|)\,  |\mathbf{u}(t)|
\\
& & + \frac{\eta}{2} \,
\left(|f(x(t),t,\eta^2t)|+|\mathbf{u}(t)|\right)\delta(|z(t)|)\\&&\times\int_{t
- 2/\eta}^{t}\int_{s}^{t}2K \, {\rm d} l \, {\rm d}s\label{ae4}
\\
& \leq & - \bar c\, \delta^2(|z(t)|)  +  K
\delta(|z(t)|)|R_\alpha(x(t),t)|\\&&+ \frac{\eta}{2}\delta(|z(t)|)
N(\eta)\, \delta(|x(t)|/2)+\, \delta(|z(t)|)\,|\mathbf{u}(t)|\\&&  +
\frac{2}{\eta} K\,
\delta(|z(t)|)\{\delta(|x(t)|/2)+|\mathbf{u}(t)|\}.\nonumber
\end{array}\]
Moreover, (\ref{m16}), the definition of $p$,  and $P_2$  give
\begin{equation}
\label{Rest}
\renewcommand{\arraystretch}{1.5}
\begin{array}{rcl}
|R_\alpha(x(t),t)|  &\leq   & \frac{\eta}{2}\int_{t -
2/\eta}^{t}\int_{s}^{t}|p(t,l)|
\, {\rm d} l\,  {\rm d}s\\
 &\leq &    \frac{2}{\eta}  \delta(|x(t)|/2) \, .
\end{array}\end{equation}
Combining these inequalities and grouping terms gives
\[\renewcommand{\arraystretch}{1.5}
\begin{array}{rcl}
\! \! \dot{V} & \leq & -\bar c\, \delta^2(|z(t)|)+
\delta(|z(t)|)\,|\mathbf{u}(t)|\\ &&+
\delta(|z(t)|)\left(\frac{4}{\eta}  K + \frac{\eta}{2}
N(\eta)\right)\{\delta(|x(t)|/2)+|\mathbf{u}(t)|\}\, .
\end{array}
\]
On the other hand, (\ref{ze1}), (\ref{Rest}), and $P_3$ give
\begin{equation}
\label{rew}
\renewcommand{\arraystretch}{1.1}
\begin{array}{rcl}
|z(t)| & \; \geq\;  &   |x(t)| - \frac{\bar{\bar c}}{\eta}  |x(t)|\;
\; \geq \; \; \frac{1}{2} |x(t)|
\end{array}
\end{equation}
when $\eta \geq \max\{2 \bar{\bar c},\eta_o\}$. Since $\delta\in
{\mathcal  K}$, this gives
\begin{equation}\label{tohere}
\renewcommand{\arraystretch}{1.1}
\begin{array}{rcl}
\dot{V} & \leq & \left(- \bar c +  \frac{4}{\eta}   K +
\frac{\eta}{2} N(\eta)\right) \delta^2(|z(t)|) \\&&+
\left(\frac{4}{\eta} K + \frac{\eta}{2}
N(\eta)+1\right)\delta(|z(t)|)|\mathbf{u}(t)| .
\end{array}\end{equation}
Setting $\chi(s)=\frac{\bar c}{4}\delta(s/2)$, it  follows from
(\ref{rew})-(\ref{tohere}) that
\begin{equation}\label{itfollows}\renewcommand{\arraystretch}{1.5}\begin{array}{l}
|\mathbf{u}|_\infty\; \le\;  \chi(|x(t)|) \; \; \; \Rightarrow\; \;
\; |\mathbf{u}|_\infty\; \le\;  \chi(2|z(t)|)\; \; \;\\
\Rightarrow\; \; \; \dot V\le \left(-\frac{3\bar
c}{4}+\frac{8}{\eta} K+\eta
N(\eta)\right)\delta^2(|z(t)|).\end{array}
\end{equation}
Setting $V^{\scriptscriptstyle
[\alpha]}(x,t):=V(x+R_\alpha(x,t),t)$, we see the derivative $\dot
V=V_t(z,t)+V_\xi(z,t)\dot z$ of $V(z,t)$ along (\ref{dis}) satisfies
\begin{equation}\label{keyformula}\dot V=V^{\scriptscriptstyle
[\alpha]}_t(x,t) + V^{\scriptscriptstyle [\alpha]}_x(x,t)\,
\{f(x,t,\alpha t)+{\mathbf u}(t)\}.\end{equation} We deduce from
(\ref{m3}), (\ref{rew}), and (\ref{itfollows}) that
 when the constant $\alpha$ (and so also  $\eta$) is sufficiently large,
\[
\renewcommand{\arraystretch}{1.7}
\begin{array}{l}
 \! \! |u|\le \chi(|x|)\; \; \; \Rightarrow\\
\! \! \! V^{\scriptscriptstyle [\alpha]}_t(x,t) +
V^{\scriptscriptstyle [\alpha]}_x(x,t)\, [f(x,t,\alpha t)+u]  \leq
 -\frac{\bar c}{2}\delta^2(|x|/2)\end{array}\]
and $\delta_1(|x|/2) \le V^{\scriptscriptstyle [\alpha]}(x,t) \le
\delta_2(|x|+2\delta(|x|/2)/\eta )$ by (\ref{Rest}).
 It follows that
$V^{\scriptscriptstyle [\alpha]}$ is an ISS Lyapunov function for
(\ref{dis}), so (\ref{dis}) is ISS for large $\alpha$, by Lemma
\ref{suff}, as claimed. The UGAS conclusion is the special case
where $\mathbf{u}\equiv 0$.  {}To prove the iISS assertion, we
instead follow the preceding argument up through (\ref{tohere})
(which is valid since that part did use the unboundedness of
$\delta$) and then substitute  $\delta(|z(t)|)|\mathbf{u}(t)| \le
\bar c\, \delta^2(|z(t)|)/2+|\mathbf{u}(t)|^2/(2\bar c)$ into
(\ref{tohere}) and bound the resulting coefficient of $|{\mathbf
u}(t)|^2$ to show that $V^{\scriptscriptstyle [\alpha]}$ is an iISS
Lyapunov function for (\ref{dis}) for sufficiently large $\alpha$
(by again using (\ref{Rest}) and (\ref{rew}) and taking
$\nu(s)=\delta^2(s/2)\bar c/4$ and $\Delta(s)=\bar r s^2$ for a
suitable $\bar r>0$), which implies that (\ref{dis}) is iISS for
large $\alpha$, by Lemma \ref{suff}.

We turn next to the special case where (\ref{m2}) is UGES.  Let $V$
satisfy the requirements of Lemma \ref{k02} above for (\ref{m2}),
 and let
$x(t)$ be any trajectory for (\ref{dis}) for any control
$\mathbf{u}\in {\mathcal  U}$ {}starting at $x(t_o)=x_o$. Define
$z(t)$ by (\ref{ze1}).
 Arguing  {}as before except with this new choice of
 $V$ shows $|R_\alpha(x(t),t)|\le   2K|x(t)|/\eta$ (by taking
 $\delta(s)=K s$ in (\ref{Rest}))
 and
  that (\ref{keyformula}) satisfies (by $P_3$ and (\ref{relate}))
\[
\renewcommand{\arraystretch}{1.1}
\begin{array}{rcl}
\dot{V} & \leq & - |z(t)|^2 +c_3|z(t)||\mathbf{u}(t)|\\&&+ c_3\,
|z(t)|\left(\frac{4}{\eta} K^2 + \frac{\bar{\bar c}\, \eta}{2}
N(\eta)\right) \{|x(t)|+|\mathbf{u}(t)|\}\\ & \leq & \left(- 1 +
\frac{8}{\eta} c_3 K^2 + \bar{\bar c}\, \eta c_3 N(\eta)\right)
|z(t)|^2\\ &&+ c_3|z(t)|\left(\frac{4}{\eta} K^2+\frac{\bar{\bar
c}\, \eta}{2}N(\eta)+1\right)|{\mathbf u}(t)|,
\end{array}\] since $|x(t)|\le 2|z(t)|$  for large $\eta$ as before.
If we now define $\tilde \chi\in {\mathcal  K}_\infty$ by $ \tilde
\chi(s) = s/\{8(1+c_3)\}$,  then we deduce as in the UGAS case that
if $\eta$ is large enough, and if $|\mathbf{u}|_\infty\le \tilde
\chi(|x(t)|)$ for all $t$, then we also have $|\mathbf{u}|_\infty\le
\tilde \chi(2|z(t)|)$ for all $t$ and $\dot{V}   \leq   - |z(t)|^2/2
   \leq   - V(z(t),t)/(2c_2)$.  This gives
$V(z(t),t)   \; \leq\;    V(z(t_o),t_o) e^{-(t-t_o)/(2c_2)}$, so
\[\frac{c_1}{4}|x(t)|^2 \; \le\; c_1|z(t)|^2 \; \leq\; V(z(t),t)  \;
\leq\;  {}c_2 |z(t_o)|^2 e^{- \frac{t-t_o}{2 c_2}}\]  so our
estimate on $|R_\alpha(x(t),t)|$  and the form of $z(t)$ give
\begin{equation}
\label{1x4a}
\renewcommand{\arraystretch}{1.1}
\begin{array}{rcl}
{}|x(t)| & \leq & \sqrt{\frac{4 c_2}{c_1}} \left(1 + \frac{2
K}{\eta}\right)|x(t_o)| e^{- \frac{t-t_o}{4 c_2}}.
\end{array}
\end{equation}
We conclude as before that if   (\ref{m2}) is UGES, then, when the
constant $\alpha>0$ is large enough, (\ref{m1}) is also UGES and
(e.g., by the proof of \cite[Lemma 2.14]{SW95}) (\ref{dis}) is ISES,
which proves our theorem.

\begin{rem}\rm
\label{gissue} The method we used in the proof of Theorem \ref{lef2}
can be used to prove the ISS property for (\ref{fastcontrol}) under
appropriate growth assumptions on the matrix-valued function
$g:{\mathbb R}^n\times {\mathbb R}\times {\mathbb R}\to{\mathbb
R}^{n\times m}$. Clearly, some growth condition on $g$ is needed and
linear growth of $g$ is not enough, since $\dot x=-x+xu$ is not ISS.
 One way to extend our theorem to (\ref{fastcontrol}) is to add the hypotheses that
{}$g$ is $C^1$ and that
 there is a
constant $c_o>1$ such that for all $ t\in {\mathbb R}$, $x\in
{\mathbb R}^n$, and  $\alpha>0$, $||g(x,t,\alpha t)||\le
c_o+(\delta(|x|/2))^{1/2}$, where $||\cdot ||$ is the $2$-norm on
${\mathbb R}^{n\times m}$ and $\delta\in {\mathcal  K}_\infty$
satisfies $P_1$-$P_3$ for some Lyapunov function $V$ for (\ref{m2}).
Applying the first part of the proof of Theorem \ref{lef2} except
with the new ${\mathcal K}_\infty$ function
\begin{equation}\label{nc}
\chi(s)=\frac{\bar c\, \delta(s/2)}{4\{c_o+\sqrt{\delta(s/2)}\}},
\end{equation} we  then conclude as before that (\ref{fastcontrol}) is ISS for
sufficiently large $\alpha>0$. If instead $\delta\in {\mathcal K}$
is bounded, then (\ref{fastcontrol}) is iISS when $\alpha$ is
sufficiently large, by our earlier argument.
\end{rem}

\begin{rem}\label{newrk} \rm

The decay requirement  (\ref{m3}) on $N\in {\mathcal  M}$ from
Theorem \ref{lef2}  can be relaxed, as follows. We assume the  flow
map $\phi$ of  (\ref{m2}) satisfies the UGES conditions
(\ref{gasdef})-(\ref{betachoice}) for some $D>1$ and $\lambda\in
(0,K)$,  where $K$ satisfies (\ref{m16}), and  we let $V$ be as in
Lemma \ref{k02}.  We can choose the
 constant $c_3$ in (\ref{m5}) to be
\begin{equation}\label{c3c}
c_3=\frac{4D(\Theta-1)}{(K-\lambda)},\; \; \; {\rm where}\; \;
\Theta=(\sqrt{2}D)^{K/\lambda-1},
\end{equation}
by the proof of Lemma \ref{k02}.  It follows from our argument above
that Theorem \ref{lef2} remains true for cases where (\ref{m2}) is
UGES  if (\ref{m3}) is relaxed to
\begin{equation} \label{relax}  \exists \eta^\star>0 \;
\,  {\rm s.t.}\; \,  \sup_{\eta\ge \eta^\star}\eta N(\eta) <
\frac{K-\lambda}{11D(\Theta-1)\bar{\bar c}}.\end{equation}
 A similar relaxation can be made in the more general UGAS setting covered by Theorem \ref{lef2}.
\end{rem}

\section{Proof of Theorem \ref{altthm}}
\label{altproof} We first set $\dot U^{\scriptscriptstyle [\alpha]}=
U^{\scriptscriptstyle [\alpha]}_t(x,t)+ U^{\scriptscriptstyle
[\alpha]}_x(x,t)f(x,t,\alpha t)$ for all $x\in {\mathbb R}^n$, $t\ge
0$, and $\alpha>0$. Using Assumptions H1.  and H3. and (\ref{fua})
with $p(t,l)$ independent of $t$ gives
\begin{equation}
\label{jb6}
\renewcommand{\arraystretch}{1.1}
\begin{array}{rcl}
\dot{U}^{[\alpha]} & \leq & - W(x,t) + p(\alpha t) \Theta(x,t) -
p(\alpha t) \Theta(x,t)\\&& + \left(\int_{t - 1}^{t} p(\alpha l)
dl\right) \Theta(x,t)
\\
& & - \left(\int_{t - 1}^{t} \left(\int_{s}^{t} p(\alpha l)
dl\right) ds\right) \left(\frac{\partial \Theta}{\partial x} f +
\frac{\partial \Theta}{\partial t}\right)
\\
& \leq & - W(x,t) + \left|\int_{t - 1}^{t} p(\alpha l) dl\right|
\frac{1}{c} W(x,t)
\\
& & + \left|\int_{t - 1}^{t} \left(\int_{s}^{t} p(\alpha l)
dl\right) ds\right| \frac{1}{c} W(x,t)
\end{array}
\end{equation}
along trajectories of (\ref{m1}), where we omit the argument
$(x,t,\alpha t)$ of $f$. For any $\alpha>0$, $t\ge 0$, and $s\in
[\alpha t-\alpha, \alpha t]$, Assumption H2. gives
\begin{equation}
\label{zeroing} \int_s^{\alpha t} p(l){\rm
d}l=\int_s^{\bar\tau(s)}p(l){\rm d}l+ \int_{\underline \tau(\alpha
t)}^{\alpha t}p(l){\rm d}l,\end{equation} where $\bar
\tau(u):=\min\{kT: k\in {\mathbb Z}, kT\ge u\}$ and $\underline
\tau(u):=\max\{kT:k\in {\mathbb Z}, kT\le u\}$. The proof of
(\ref{zeroing}) uses the facts that both endpoints of the interval
$[\min\{\bar\tau(s),\underline \tau(\alpha t)\},
\max\{\bar\tau(s),\underline \tau(\alpha t)\}]$ are integer
multiples of $T$ and that, according to Assumption H2, the integral
of $p$ over any interval whose endpoints are integer multiples of
$T$ is equal to zero. Choosing $p_{\rm max}$ to be any global bound
on $|p(l)|$ over all $l\in {\mathbb R}$, (\ref{zeroing})  gives
\begin{equation}
\label{bounda} \left\vert\int_s^{\alpha t} p(l){\rm d}l\right\vert\;
\le \; 2Tp_{\rm max}\; \; \; \; \forall s\in [\alpha t-\alpha,
\alpha t]
\end{equation}
for all $t\ge 0$ and $\alpha>0$.  Hence,   for all $t\ge 0$,
\begin{eqnarray}
\label{termest} \left\vert\int_{t - 1}^{t} p(\alpha l) dl\right\vert
& = & \frac{1}{\alpha}\left\vert\int_{\alpha t - \alpha}^{\alpha t}
p(l)
dl\right\vert\le \frac{2Tp_{\max}}{\alpha}\nonumber\\
\left|\int_{t - 1}^{t} \left(\int_{s}^{t} p(\alpha l) dl\right)
ds\right| & = & \frac{1}{\alpha^2} \left|\int_{\alpha t -
\alpha}^{\alpha t} \left(\int_{s}^{\alpha t} p(l) dl\right)
ds\right|\nonumber
\\
& \leq & \frac{1}{\alpha} \sup_{s \in [\alpha t - \alpha, \alpha
t]}\left|\int_{s}^{\alpha t} p(l) dl\right|\nonumber\\ &\le &
2Tp_{\max}/\alpha
\end{eqnarray}
  Using these estimates, (\ref{jb6}) gives $\dot
U^{\scriptscriptstyle [\alpha]}\le -W(x,t)/2$ for all $x\in {\mathbb
R}^n$ and $t\ge 0$, as long as $\alpha>8Tp_{\max}/c$. Since $W\in
{\rm UPD}$, this gives the Lyapunov decay estimate.  By Assumption
H3. and (\ref{termest}),
 $U^{\scriptscriptstyle [\alpha]}\in {\rm UPPD}$ for large
enough constants $\alpha>0$.  This proves the theorem.

\section{Illustrations of Theorem \ref{lef2}}

\label{illu}  We next illustrate how Theorem \ref{lef2} extends the
 results of \cite{K02} and \cite{PA02}. In \cite{PA02}, the
limiting dynamics (\ref{m2}) are assumed to be {}uniformly locally
exponentially stable. {}In our first example, {}(\ref{m2}) is UGAS
but not necessarily UGES. We next consider a class of systems
(\ref{m1}) from \cite{PA02} {}from identification where the limiting
dynamics (\ref{m2}) is linear and exponentially stable. For these
systems, our work {}complements \cite{PA02} by providing formulas
for Lyapunov functions for (\ref{m1}) that are expressed in terms of
the quadratic Lyapunov functions for the limiting dynamics and
 that have the additional desirable
property that they are also ISS Lyapunov functions for
(\ref{fastcontrol}) for suitable functions $g$.  Finally, we  apply
our results to a friction model for a mass-spring dynamics. In all
three examples, the limiting dynamics has a simple  Lyapunov
function structure so our results give explicit Lyapunov functions
for the original rapidly time-varying dynamics.

\subsection{Application to a UGAS dynamics that is not UGES}\label{exa1}
Consider the following variant of the scalar example on
\cite[p.53]{PA02}:
\begin{equation}\label{nonUGES}
\begin{array}{l}
\dot x =f(x,t,\alpha t)= \\
-\sigma_1(x)[2+\sin(t+\cos(\sigma_2(x)))]\{1+10\sin(\alpha
t)\}\end{array}
\end{equation}
where $\sigma_1,\sigma_2:{\mathbb R}\to {\mathbb R}$ are $C^1$
functions such that $\sigma_1$ is odd,
$\sup\{|\sigma'_1(x)|+|\sigma_1(x)\sigma'_2(x)|: x\in {\mathbb
R}\}<\infty$, $\sigma_1\in {\mathcal  K}$ on $[0,\infty)$, and
$\sigma''_1(s)\le 0$ for all $s>0$. One easily verifies the
hypotheses of Theorem \ref{lef2} using
\[\begin{array}{l}
\bar{f}(x,t):=-\sigma_1(x)[2+\sin(t+\cos(\sigma_2(x)))],\\
V(x,t)\, \equiv\, \bar V(x)\, :=\,  \int_0^{x}\sigma_1(s){\rm
d}s,\\
\delta(s):=33\sigma_1(2s),\; \; {\rm and}\; \; N(\eta):=60/\eta^2 \;
{\rm for\ large\ }  \eta.\end{array}\] This allows e.g.
$\sigma_1(s)=\sigma_2(s)=\arctan(s)$ in which case (\ref{m2}) is
UGAS but not UGES because $|\dot x(t)|\le 2\pi$ along all of its
trajectories $x(t)$. Condition $P_1$ follows because
$\sigma_1(2s)\le 2\sigma_1(s)$ for all $s\geq 0$, which holds
because $\sigma''_1(s)\leq 0$ for all $s\ge 0$. Theorem \ref{cor1}
then gives the following iISS Lyapunov function for (\ref{dis}) for
large $\alpha>0$:
\begin{equation}\begin{array}{l}\! \! \! \! \! \! \!
\bar V\left(\xi+5\sqrt{\alpha}\, \sigma_1(\xi) \int_{
t-\frac{2}{\sqrt{\alpha}}}^t\int_s^t \mu(\xi,l)\sin(\alpha l) \,
{\rm d} l\, {\rm d}s\right),\end{array}
\end{equation}
where $\mu(\xi,l):=2+\sin(l+\cos(\sigma_2(\xi)))$.  In particular,
this is a Lyapunov function for $\dot x=f(x,t,\alpha t)$, and it is
also an ISS Lyapunov function for (\ref{dis}) if $\delta\in
{\mathcal K}_\infty$ (e.g. if $\sigma_1(s)={\rm sgn}(s)\ln(1+|s|)$
for $|s|\ge 1$ and $\sigma_2(s)=\arctan(s)$).
 Our conditions on the $\sigma_i$'s  cannot be
omitted even if the limiting dynamics {}(\ref{m2}) is UGES; see
\cite[$\S$8.2]{PA02}. For example, if $\sigma_1(x) = x$ and
$\sigma_2(x)=x^2$, then (\ref{m2}) is UGES, but (\ref{nonUGES}) is
only shown to be {\em locally} exponentially stable for large
$\alpha>0$;  see \cite{PA02}. This does not contradict our theorem
because in that case (\ref{m16}) would be violated.

\subsection{A system arising in identification}
Consider the following variant of the example in \cite[Section
8.1]{PA02}:
\begin{equation}
\label{pertid} \dot x\; \; =\; \; f(\alpha t)\, m(t)\, m^T(t)\,
x+g(x,t,\alpha t)\, u,\end{equation} with state $x\in {\mathbb R}^n$
and inputs $u\in {\mathbb R}^m$, where we assume
\begin{itemize}\item[]\begin{itemize}
\item[i.] $f:{\mathbb R}\to {\mathbb R}$ is  bounded and continuous
and admits  a $o(s)$ function $M$  and a constant $f^\star<0$ for
which
\begin{enumerate}
\item $ f^\star=\lim_{T\to +\infty}\frac{1}{2T}\int_{-T}^Tf(s){\rm
d}s$
\item $\vert\int_{t_1}^{t_2}[f(s)-f^\star]{\rm d}s\vert\le M(t_2-t_1)$
if $t_2\ge t_1$\smallskip
\end{enumerate}
\item[ii.]$m:{\mathbb R}\to {\mathbb R}^n$ is continuous and
satisfies $|m(t)|=1$ for all $t\in {\mathbb R}$, and there exist
constants $\alpha', \beta', \tilde c>0$ such that for all $t\in
{\mathbb R}$, $\alpha>0$, and $x\in {\mathbb R}^n$, we have:
\begin{enumerate}
\item
$\alpha' I\le \int_t^{t+\tilde c}m(\tau)m^T(\tau){\rm d}\tau\le
\beta' I$
\item $||g(x,t,\alpha t)||\le  \beta'\{1+\sqrt{|x|}\}$\smallskip
\end{enumerate}
\item[iii.] $g:{\mathbb R}^n\times {\mathbb R}\times  {\mathbb R}\to
{\mathbb R}^{n\times  m}$ is continuous and is $C^1$ in $x$, and
there exists a constant $K>1$ such that $|\partial g_{ij}(x,t,\alpha
t)/\partial x|\le K$ $\forall x\in {\mathbb R}^n$, $t\ge 0$, and
$\alpha>0$, and each component $g_{ij}$ of $g$.
\end{itemize}\end{itemize} Here $I$ denotes the $n\times n$ identity
matrix. Notice that we allow $f$ to take both positive and negative
values. The special case of (\ref{pertid}) where $g\equiv 0$ was
studied in \cite{PA02} where it is shown that the corresponding
rapidly varying dynamics $\dot x= f(\alpha t) m(t) m^T(t) x$
satisfies the hypotheses of Theorem \ref{lef2}  with the UGES
dynamics
\begin{equation}
\label{limid} \dot x=\bar f(x,t):=f^\star m(t)\, m^T(t)\,
x\end{equation} and with $\delta$ of the form $\delta(s)=\bar r s$
for a constant $\bar r>0$.
  The particular case
of (\ref{pertid}) in which $\dot x= -m(t)\, m^T(t)\, x$ has been
extensively studied in systems identification; see \cite{PA02}.
However, these earlier results do not provide explicit ISS Lyapunov
functions for (\ref{pertid}).  On the other hand, the following
lemma provides an explicit Lyapunov function for (\ref{limid}):

\begin{lem}\label{plem}
Let assumptions i.-iii. hold. If we choose
\begin{equation}
\label{pchoice} P(t)=\kappa I+\int_{t-\tilde c}^t\int_s^t
m(l)m^T(l)\, {\rm d}l\, {\rm d}s,\end{equation} where $\kappa=
\tilde c/(2 |f^\star|) + \frac{1}{4 \alpha'} \tilde c^4 |f^\star|$
then $V(x,t)=x^TP(t)x$ is a Lyapunov function for (\ref{limid}) for
which $2V/\alpha'$ satisfies the requirements of Lemma
\ref{k02}.\end{lem}

To prove Lemma \ref{plem}, we apply (\ref{fua})-(\ref{fub}) with
$\tau=\tilde c$ and  $p(t,l)\equiv m(l)m^T(l)$ and   group terms to
check that the derivative of $V$ along trajectories of (\ref{limid})
satisfies
\[\begin{array}{rcl}
\dot{V} & = & (2 f^\star \kappa + \tilde c) x^T m(t) m^T(t) x
\\&&+ 2 f^\star x^T \left[\int_{t-\tilde c}^t\int_s^t
m(l)m^T(l)\, {\rm d}l\, {\rm d}s\right]  m(t) m^T(t) x
\\
& & - x^T \left[\int_{t-\tilde c}^t m(l)m^T(l)\, {\rm d}l \right]
x\\ \dot{V} & \leq & (2 f^\star \kappa + \tilde c) |m^T(t) x|^2-
\alpha' |x|^2
\\&& + 2 |f^\star| |x| \left[\int_{t-\tilde
c}^t\int_s^t |m(l)|^2 {\rm d}l {\rm d}s\right]  |m(t)| |m^T(t) x|
\\
& \leq & \! (2 f^\star \kappa + \tilde c) |m^T(t) x|^2 + \tilde c^2
|f^\star| \, |x| |m^T(t) x| - \alpha' |x|^2
\end{array}
\]
everywhere, by ii.(1). By the triangle inequality, \[  \tilde c^2
|f^\star| |x| |m^T(t) x| \leq \frac{1}{2}\alpha' |x|^2 + \frac{1}{2
\alpha'} \tilde c^4 |f^\star|^2 |m^T(t) x|^2
\]
so  $\dot{V}  \leq  \left(2 f^\star \kappa + \tilde c + \frac{1}{2
\alpha'} \tilde c^4 |f^\star|^2\right) |m^T(t) x|^2 -
\frac{1}{2}\alpha' |x|^2$ holds everywhere.  Recalling that
$|m(t)|=1$ everywhere and that $P$ is everywhere positive definite,
the fact that $2V/\alpha'$ satisfies the requirements of Lemma
\ref{k02} follows immediately from our choice of
 $\kappa$.
 Remark \ref{gissue} now gives:
\begin{cor} Let (\ref{pertid}) satisfy conditions i.-iii. above and
let $V$ be as in Lemma \ref{plem}.  Then there exists a constant
$\alpha_o>0$ such that for each constant $\alpha>\alpha_o$,
\[
V^{[\alpha]}(x,t):=V\left(\left[I-\frac{\sqrt{\alpha}}{2}\int_{t-2/\sqrt{\alpha}}^t\int_s^t
D(l)\, {\rm d}l\, {\rm d}s\right]x,t\right)
\]
where $D(l)=(f(\alpha l)-f^\star)m(l)m^T(l)$ is an ISS Lyapunov
function for (\ref{pertid}).\end{cor}

\subsection{Friction example}

The following one degree-of-freedom mass-spring system from
\cite{DDNZ00} arises in the control of mechanical systems in the
presence of friction:
\begin{equation}
\label{frictionexample}
\begin{array}{rcl}
\dot{x}_{1} &=&x_{2} \\
\dot{x}_{2} &=&-\sigma _{1}(\alpha t)x_{2}-k(t)x_{1}+u\\&&-\left\{
\sigma _{2}(\alpha t)+\sigma _{3}(\alpha t)e^{-\beta _{1}\mu
(x_{2})}\right\} {\rm sat} (x_{2})\end{array}
\end{equation}
where $x_{1}$ and $x_{2}$ are the mass position and velocity,
respectively; $\sigma _{i}$, $i=1,2,3$ denote positive time-varying
viscous, Coulomb, and static friction-related coefficients,
respectively; $\beta _{1}$ is a positive constant corresponding to
the Stribeck effect; $\mu (\cdot )$ is a positive definite function
also related to the Stribeck effect; $k$ denotes a positive
time-varying spring stiffness-related coefficient; and sat$(\cdot )$
denotes any continuous function having these properties:
\begin{equation}\label{sats}\renewcommand{\arraystretch}{1.2}
\begin{array}{l} \! \! \! \! \!
\! {\rm (a)}\; \; {\rm sat}(0)=0, \ \ \  {\rm (b)}\; \;  \xi \, {\rm
sat}(\xi )\geq 0\; \; \forall \xi\in {\mathbb R},\\ \! \! \! \! \!
\! {\rm (c)}\; \; \lim\limits_{\xi \rightarrow + \infty }{\rm
sat}(\xi )=+ 1,  \; \; \; {\rm (d)}\; \; \lim\limits_{\xi
\rightarrow - \infty }{\rm sat}(\xi )=- 1
\end{array}\end{equation} We model the saturation as the
differentiable function
\begin{equation}
{\rm sat}(x_{2})=\tanh (\beta _{2}x_{2}),  \label{sat}
\end{equation}
where $\beta _{2}$ is a large positive constant. Note for later use
that  $|{\rm sat}(x_{2})|\le \beta_2|x_2|$ for all $x_2\in {\mathbb
R}$. We assume the friction coefficients vary in time faster than
the spring stiffness coefficient so we restrict to cases where
$\alpha>1$.

Our precise mathematical assumptions on  (\ref{frictionexample})
are: $k$ and the $\sigma_i$'s are bounded $C^1$ functions;   $\mu $
has a globally bounded derivative; and there exist  constants
$\tilde \sigma_i$, with $\tilde \sigma_1>0$ and $\tilde \sigma_i\ge
0$ for $i=2,3$, and a $o(s)$ function $s\mapsto M(s)$ such that
\begin{equation}
\label{sigma1} \left| \int_{t_{1}}^{t_{2}}\left( \sigma
_{i}(t)-\tilde{\sigma}_i\right) \mathrm{d}t\right| \ \le \
M(t_{2}-t_{1}),\; \; \; i=1,2,3
\end{equation}
for all $t_{1},t_{2}\in {\mathbb R}$ satisfying $t_{2}>t_{1}$.
Although the $\sigma_i$'s are positive for physical reasons, we will
not require their positivity in the sequel.
We show (\ref{frictionexample}) satisfies the requirements of the
version of Theorem \ref{lef2} from Remark \ref{gissue} (with
$\delta(s)=\bar r s$ for a constant $\bar r$) when (\ref{m2}) is
\begin{eqnarray} \label{fbf} \, \, \dot{x}_{1} &=&x_{2}\\ \, \,
\dot{x}_{2} &=&-\tilde{\sigma}_1x_{2}- \left\{\tilde \sigma_2+\tilde
\sigma _{3}e^{-\beta _{1}\mu (x_{2})}\right\} {\rm sat}
(x_{2})-k(t)x_{1}\, ,\nonumber\end{eqnarray} assuming this
additional condition whose physical interpretation is that the
spring stiffness is nonincreasing:
\[
\exists k_o, \bar k>0\; \; {\rm s.t.}\; \; k_o \le k(t)\le \bar k\;
\; {\rm and}\; \; k'(t)\le 0 \; \; \forall t\ge 0.
\]
To this end, set $S:=\tilde
\sigma_1+(\tilde\sigma_2+\tilde\sigma_3)\beta_2$ and
\begin{equation}
\label{qu1} V(x,t) = A (k(t)x_1^2 + x_2^2) + x_1 x_2,
\end{equation}
where $A:=1+1/k_o+[1+S^2/k_o]/\tilde \sigma_1$. Since  $A\bar k\ge
1$, $\frac{1}{2}(x^2_1+x^2_2)\le V(x,t)\le A^2 \bar
k(|x_1|+|x_2|)^2$ for all $x \in {\mathbb R}^2$ and $t\ge 0$. Also,
since $k'\le 0$ everywhere, the derivative $\dot
V=V_t(x,t)+V_x(x,t)\bar f(x,t)$  along trajectories of (\ref{fbf})
satisfies
\[\renewcommand{\arraystretch}{1.25}
\begin{array}{l}
\! \! \! \dot V\le V_x(x,t)\bar f(x,t)
=[2Ak(t)x_1+x_2]x_2-[2Ax_2+x_1]\\\; \; \; \; \; \; \;
\times\{\tilde{\sigma}_1x_{2}+ \left[\tilde \sigma_2+\tilde \sigma
_{3}e^{-\beta _{1}\mu (x_{2})}\right] {\rm sat} (x_{2})+k(t)x_{1}\}
\end{array}\]
and therefore, by grouping terms, we also have
\[
\renewcommand{\arraystretch}{1.5}
\begin{array}{l}
\dot V\le -k_0x^2_1-(2A\tilde \sigma_1-1)x^2_2+S|x_1x_2|\; \; \; \; \; {\rm (by\  (\ref{sats})(b))}\\
\le
-b|x|^2-\left[\frac{k_o}{2}x^2_1+(A\tilde\sigma_1-1/2)x^2_2-S|x_1x_2|\right]\\
= -b|x|^2-\frac{k_o}{2}\left(|x_1|-\frac{S}{k_o}|x_2|\right)^2\!
+\left(\frac{S^2}{2k_o}+\frac{1}{2}
-A\tilde \sigma_1\right)x^2_2\\
\le -b|x|^2,\;  {\rm where}\;  b:=\min\{k_o/2,A\tilde
\sigma_1-1/2\}.
\end{array}\]
The preceding inequalities imply that $V/b$ is a Lyapunov function
for (\ref{fbf}) satisfying the requirements of Lemma \ref{k02}. The
integral bound requirement {}(\ref{relate}) from our theorem follows
from (\ref{sigma1}) and the sublinear growth of $\tanh$, since the
integral bound can be verified term by term. The remaining bounds
from (\ref{m16}) follow because $\mu$ and sat have globally bounded
derivatives. We conclude that for sufficiently large constants
$\alpha>0$, (\ref{frictionexample}) admits the ISS Lyapunov function
\[
\begin{array}{l}
\! \! \! V^{\scriptscriptstyle [\alpha]}(\xi,t)=V\left(\xi_1, \xi_2
+ \frac{\sqrt{\alpha}}{2}\int_{t-\frac{2}{\sqrt{\alpha}}}^t\int_s^t
\Gamma_\alpha(l,\xi)\, {\rm d} l\, {\rm d}s,t\right)
\end{array}
\]
where $V$ is the Lyapunov function (\ref{qu1}) for (\ref{fbf}) and
\begin{eqnarray}
\Gamma_\alpha(l,\xi):=\{\sigma_1(\alpha
l)-\tilde\sigma_1\}\xi_2+\mu_\alpha(l,\xi)\tanh(\beta_2 \xi_2)\\
\mu_\alpha(l,\xi):=\sigma_2(\alpha l)-\tilde
\sigma_2+(\sigma_3(\alpha l)-\tilde
\sigma_3)e^{-\beta_1\mu(\xi_2)}\end{eqnarray} so
(\ref{frictionexample}) is ISS for large enough $\alpha>0$, by
Remark \ref{gissue}.

\begin{rem}\rm
\label{const}  The preceding construction simplifies considerably if
$\sigma_2$ and $\sigma_3$ in (\ref{frictionexample}) are positive
constants. In that case, the limiting dynamics (\ref{m2}) can be
taken to be
\begin{eqnarray*}
\dot{x}_{1} \; &\; =\; &\; x_{2} \\
\dot{x}_{2} \; &\; =\; &\; -\tilde{\sigma}_1x_{2}-\left\{\sigma_2+
\sigma _{3}e^{-\beta _{1}\mu (x_{2})}\right\} {\rm sat} (x_{2})
-k(t)x_{1}
\end{eqnarray*}
and  the  ISS Lyapunov function for (\ref{frictionexample}) becomes
\[
\begin{array}{l}
\! \! \! V\left(\xi_1, \xi_2\left(1+\frac{\sqrt{\alpha}}{2}\int_{
t-\frac{2}{\sqrt{\alpha}}}^t\int_s^t \{\sigma_1(\alpha l)-\tilde
\sigma_1\}\, {\rm d} l\, {\rm d}s\right) ,t\right)
\end{array}
\]
with $V$ defined by  (\ref{qu1}), since the
$\mu_\alpha(l,\xi)\tanh(\beta_2 \xi_2)$ term in the difference
$f-\bar f$ is no longer present.
\end{rem}

\section{Illustrations of Theorem \ref{altthm}}

 Simple calculations show that Theorem \ref{altthm} applies to
 (\ref{nonUGES}) and (\ref{pertid}) without controls (with the same choices of $V$
 that we used in our earlier discussions of  those systems), assuming in
 the latter case that requirements i.-iii.
hold and for instance $f$ is a suitable periodic function and
$\dot{m}(t)$ is bounded.  We next show how Theorem \ref{altthm} also
applies to cases that are not tractable by Theorem \ref{lef2}.

\label{altillus}
\subsection{Dynamics that are not globally Lipschitz} Simple
calculations that we omit because of space constraints show that the
one-dimensional dynamics
\begin{equation}
\label{ngs} \dot{x} \; = \; - x^3 + 10 \cos(\alpha t) \frac{x^3}{1 +
x^2}
\end{equation}
satisfy H. using $V(x,t)\equiv x^4/4$, $W(x,t)\equiv x^6$,
$\Theta(x,t)\equiv x^6/(1+x^2)$, $p(t)=10\cos( t)$, $T=2\pi$, and
small enough $c>0$, so (\ref{ngs}) has the global Lyapunov function
\[ U^{\scriptscriptstyle [\alpha]}(x,t)=
\frac{x^4}{4} - \left(\displaystyle\int_{t - 1}^{t}
\left(\displaystyle\int_{s}^{t} 10 \cos(\alpha l) dl\right)
ds\right)\frac{x^6}{1 + x^2}\] when $\alpha>0$ is large enough.
However, (\ref{ngs}) is not covered by  Theorem \ref{lef2} since it
is not globally Lipschitz in $x$.

\subsection{Dynamics with unknown functional parameters} Consider
the nonautonomous scalar system
\begin{equation}
\dot{x}=p(\alpha t)\frac{x^{2}}{1+x^{2}}+u,  \label{ex}
\end{equation}
where $p$ is an unknown, fast time-varying parameter satisfying
$\left| p(l)\right| \leq a_{m}$ for all $l$ and some constant
$a_m>0$ and admitting a constant $T>0$ such that  Assumption H2.
holds for all $k\in {\mathbb Z}$.  We assume  the control  $u$ is
amplitude limited in the sense that $\left| u\right| \leq u_{m}$ for
some constant $u_m>0$. We  show that the saturated state feedback
\begin{equation}
u=-u_{m}\arctan (x)\quad   \label{u}
\end{equation}
renders (\ref{ex}) UGAS. (A similar argument shows  that
$u=-u_{m}\arctan (Rx)$ stabilizes (\ref{ex}) for any constant
$R>1$.) For simplicity, let $a_{m}=10$ and $u_{m}=2$.  The
derivative of $V(x)=x^2/2$ along (\ref{ex}) in closed loop with
(\ref{u}) is
\begin{equation}
\dot{V}=-2x\arctan (x)+p(\alpha t)\frac{x^{3}}{1+x^{2}}.
\label{Vdot}
\end{equation}
Simple calculations allow us to verify the hypotheses of Theorem
\ref{altthm} for the closed loop system using  $ W(x,t)\equiv
2x\arctan (x)$, $\Theta(x,t)\equiv x^3/(1+x^{2})$, and small enough
$c>0$, so we know (\ref{u}) indeed uniformly globally asymptotically
stabilizes (\ref{ex}) with control Lyapunov function
\begin{equation}
U^{\scriptscriptstyle [\alpha]}(x,t)=\dfrac{1}{2}x^{2}-\left(
\int_{t-1}^{t}\int_{s}^{t}p(\alpha l )\, {\rm d}l\,  {\rm d}s\right)
\frac{x^{3}}{1+x^{2}}  \label{clf}
\end{equation}
when $\alpha>0$ is sufficiently large.

\section{Conclusions}
\label{concl} {}Hypotheses similar to those  of \cite[Theorem
3]{PA02} are sufficient for uniform global  exponential stability of
rapidly time-varying nonlinear systems. We provided complementary
results to \cite{PA02} by establishing uniform global asymptotic
 {}stability of fast time-varying
dynamics without requiring {}local exponential stability of
 the limiting dynamics, and by
constructing global Lyapunov functions for  fast time-varying
systems. Our Lyapunov constructions are new even in the special case
{}of exponentially stable dynamics, and are input-to-state stable
Lyapunov functions when the dynamics are control affine, under
appropriate conditions.  Our results apply to dynamics that are not
necessarily uniformly Lipschitz in the state. We illustrated our
methods using a friction example.

\thebibliography{xx}

\harvarditem[Angeli {\em et al.}]{Angeli {\em et al.}}{2000}{ASW00}
Angeli, D., Sontag, E.D., \& Wang, Y. (2000). A characterization of
integral input to state stability. {\it IEEE Transactions  on
Automatic Control}, 45(6),  1082-1097.

\harvarditem[Bacciotti \& Rosier]{Bacciotti \& Rosier}{2005}{BR05}
 Bacciotti, A.,  \& Rosier, L. (2005).  {\it Liapunov Functions and Stability
 in Control Theory} (2nd ed.). New York: Springer Verlag.

\harvarditem[de Queiroz {\em et al.}]{de Queiroz {\em et al.}}{2000}
{DDNZ00} de Queiroz, M.S., Dawson, D.M.,  Nagarkatti, S.,  \& Zhang,
F. (2000). {\it Lyapunov-Based Control of Mechanical Systems.}
Cambridge, MA: Birkh\"{a}user.

\harvarditem[Edwards {\em et al.}]{Edwards, Lin \&
Wang}{2000}{ELW00} Edwards, H., Lin, Y., \& Wang, Y. (2000). On
input-to-state stability for time-varying nonlinear systems. {\em
Proceedings of the 39th IEEE Conference on Decision and Control},
Sydney, Australia, (pp. 3501-3506).

\harvarditem[Gr\"une {\em et al.}]{Gr\"une, Sontag \&
Wang}{1999}{GSW99} Gr\"une, L.,  Sontag, E.D., \& Wirth, F. (1999).
Asymptotic stability equals exponential stability, and ISS equals
finite energy gain - if you twist your eyes. {\em Systems and
Control Letters}, 38(2),  127-134.

\harvarditem[Hale]{Hale}{1980}{H80} Hale, J. K. (1980).\  {\it
Ordinary Differential Equations}. Malabor, FL: Krieger.

\harvarditem[Khalil]{Khalil}{2002}{K02} Khalil, H. (2002). {\it
Nonlinear Systems} (3rd ed.). Englewood Cliffs, NJ: Prentice Hall.


\harvarditem[Malisoff \& Mazenc]{Malisoff \& Mazenc}{2005}{MM05}
Malisoff, M.,  \& Mazenc, F. (2005). Further remarks on strict
input-to-state stable Lyapunov functions for time-varying systems.
{\it Automatica}, 41(11),  1973-1978.

\harvarditem[Mazenc \& Bowong]{Mazenc \& Bowong}{2004}{MB04} Mazenc,
F.,  \& Bowong, S. (2004). Backstepping with bounded feedbacks for
time-varying systems. {\it SIAM Journal on Control and
Optimization}, 43(3), 856-871.

\harvarditem[Peuteman \& Aeyels]{Peuteman \& Aeyles}{2002}{PA02}
Peuteman, J.,  \& Aeyels, D. (2002). Exponential stability of
nonlinear time-varying differential equations and partial averaging.
{\it Mathematics of Control, Signals, and Systems}, 15(1), 42-70.

\harvarditem[Sontag]{Sontag}{1989}{S89}  Sontag, E.D. (1989). Smooth
stabilization implies coprime factorization. {\it IEEE Transactions
on  Automatic Control}, 34(4), 435-443.

\harvarditem[Sontag]{Sontag}{1998}{S98} Sontag, E.D. (1998).
Comments on integral variants of ISS. {\it Systems and Control
Letters}, 34(1-2),  93-100.

\harvarditem[Sontag \& Wang]{Sontag \& Wang}{1995} {SW95} Sontag,
E.D., \& Wang, Y. (1995).
 On characterizations of the input-to-state
stability property. {\it Systems and Control Letters}, 24(5),
351-359.


\endthebibliography

\end{document}